\documentclass{amsart}
\usepackage[symbol]{footmisc}
\usepackage[T1]{fontenc}
\usepackage{amsmath,amssymb,amsthm}
\usepackage{algorithm}
\usepackage{algpseudocode}
\usepackage{graphicx}
\usepackage{subfig}
\usepackage{array}
\usepackage{tikz}
\usepackage[colorlinks=true,
            linkcolor=blue,
            citecolor=blue,
            urlcolor=blue]{hyperref}

\newtheorem{theorem}[subsection]{Theorem}
\newtheorem{lemma}[subsection]{Lemma}
\newtheorem{proposition}[subsection]{Proposition}

\theoremstyle{definition}
\newtheorem{definition}[subsection]{Definition}
\newtheorem{remark}[subsection]{Remark}

\title{Solving Admissibility for the Spatial X-Ray Transform On the Two Element Field}
\author{Mihika Dusad}
\thanks{Massachusetts Institute of Technology}

\date{December 1, 2025}

\begin{document}
\maketitle
\vspace{-24pt}

\begin{abstract}
The admissibility problem in integral geometry asks for which collections of affine subspaces the Radon transform remains injective. In the discrete setting, this becomes a purely combinatorial question about recovering a function on a finite vector space from its sums over a prescribed family of affine subspaces. In this paper, we study the spatial X-ray transform (line transform) over the finite vector spaces $\mathbb{Z}_{2}^{n}$ and give a complete structural and enumerative description of admissible line complexes in $\mathbb{Z}_{2}^{4}$.

We prove that any admissible line complex in $\mathbb{Z}_{2}^{4}$ can be obtained by taking a disjoint union of one or more odd cycles and attaching trees to the cycle vertices. Using this structural description, we carry out a systematic case-by-case enumeration of all admissible complexes in $\mathbb{Z}_{2}^{4}$ and derive an exact total count. We then generalize our approach to an algorithm that applies to $\mathbb{Z}_{2}^{n}$ for arbitrary $n$, and we then implement it to obtain the total number of admissible complexes in $\mathbb{Z}_{2}^{5}$. Our results extend previous small-dimensional classifications and provide an algorithmic framework for studying admissibility in higher dimensions. Beyond their intrinsic combinatorial interest, these structures model discrete sampling schemes for tomographic imaging, and they suggest further connections between admissibility, incidence matrices, and spectral properties of the associated graphs.
\end{abstract}

\section{Introduction}

Classically, the Radon transform expresses a function in terms of its integrals over affine subspaces and lies at the heart of integral geometry and modern tomography. It first appeared in Radon's original work on determining a function from its hyperplane integrals \cite{fundamentalRadon} and was subsequently developed in the general geometric setting of homogeneous and symmetric spaces by Helgason \cite{helgasonRadon} and others. Integral-geometric ideas can be traced further back to work of Chern on Klein spaces and integral geometry in curved settings \cite{sschern1942}. In imaging, the same transform underlies the mathematical theory of computed tomography as developed, for instance, in Natterer's monograph \cite{nattererCT}. 

In parallel with the continuous theory, one may consider Radon and X-ray-type transforms over finite fields. These were studied in particular by Koh, who obtained sharp $L^{p}$–$L^{r}$ estimates for the Radon and X-ray transforms in finite fields \cite{koh2012sharp}. In this finite-field discrete setting, the admissibility problem becomes purely combinatorial: one seeks to find for which collections of affine subspaces the discrete transform is injective.

\subsection{The discrete admissibility problem}

Let $q$ be a prime power and let $\mathbb{F}_{q}^{n}$ denote the $n$-dimensional vector space over the finite field $\mathbb{F}_{q}$. For a function $f \colon \mathbb{F}_{q}^{n} \to \mathbb{C}$ and a $k$-dimensional affine subspace $H \subset \mathbb{F}_{q}^{n}$, the discrete $k$-plane Radon transform is defined by
\[
R_{k}f(H) \;=\; \sum_{x \in H} f(x).
\]
Given a subcollection $\mathcal{C}$ of affine $k$-planes, one may restrict $R_{k}$ to $\mathcal{C}$ and ask whether $f$ is uniquely determined by the values $\{R_{k}f(H)\}_{H \in \mathcal{C}}$. If this is the case for all functions $f$, we say that $\mathcal{C}$ is \emph{admissible}. This terminology is traced back to Gelfand and collaborators, who formulated an admissibility problem for integral geometry in a variety of continuous settings; see, for example, \cite{gelfand1966} and the discussion in \cite{grinbergPHD}.

In the discrete setting, admissibility is equivalent to the invertibility of a certain linear map: the incidence matrix between points of $\mathbb{F}_{q}^{n}$ and elements of $\mathcal{C}$. Thus the problem is naturally combinatorial. When $q=2$ and $k=1$, the geometry of $\mathbb{Z}_{2}^{n}$ is particularly amenable to explicit description, and the spatial X-ray transform reduces to summation over lines.

\subsection{Admissible line complexes over $\mathbb{Z}_{2}^{n}$}

Throughout this paper, we focus on the case where $q=2$ and $k=1$. We view the finite vector space $\mathbb{Z}_{2}^{n}$ as a set of $2^{n}$ vertices, and we identify an affine line by a pair of distinct vertices along with the unique affine line passing through them. A \emph{line complex} is a subset of the set of all lines in $\mathbb{Z}_{2}^{n}$.

Admissibility for line complexes in $\mathbb{Z}_{2}^{3}$ has been studied earlier, both combinatorially and computationally \cite{Grinberg_2012}. In particular, the total number of admissible line complexes in $\mathbb{Z}_{2}^{3}$ is known to be $937440$. This classification makes heavy use of the geometry of the cube and the structure of cycles in the corresponding graph, and it is illustrated in depth in the ``scrapbook'' of inadmissible complexes compiled by Grinberg and Orhon \cite{grinberg2021}.

\begin{figure}[h]
\centering
\includegraphics[scale = 0.35]{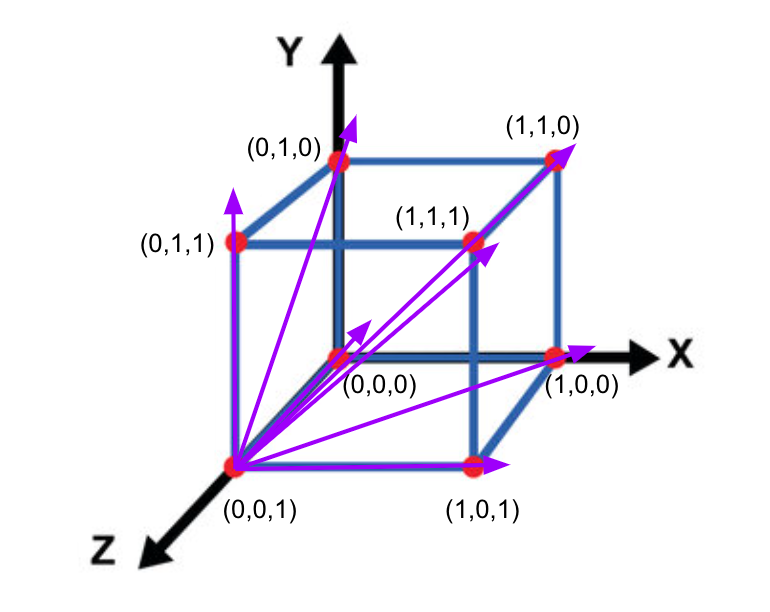}
  \caption{The Radon transform on the two element field $\mathbb{Z}_{2}^{3}$.}
  \label{fig:cube}
\end{figure}

Extending these methods directly to $\mathbb{Z}_{2}^{4}$ is not feasible via brute force. There are
\[
\binom{120}{16} \approx 3.10 \times 10^{19}
\]
possible line complexes of size $2^{4}=16$, and only a small proportion of these are admissible. Our goal is to obtain a structural characterization of admissible complexes that allows us to count them without enumerating all possibilities.

\subsection{Main results}

The first main result of this paper is a structural description of admissible line complexes in $\mathbb{Z}_{2}^{n}$ in terms of cycles and trees. We introduce a graph-theoretic viewpoint: points of $\mathbb{Z}_{2}^{n}$ are vertices, lines are edges, and a line complex is a spanning subgraph. We show that basic obstructions to admissibility arise from omitted points, isolated trees, and even cycles. These obstructions already appear in the three-dimensional analysis of \cite{Grinberg_2012,grinberg2021}, where trees and even cycles are used to generate explicit inadmissible examples. We prove the following theorem.

\begin{theorem}\label{thm:structure-overview}
Let $\mathcal{C}$ be a line complex in $\mathbb{Z}_{2}^{n}$. Suppose that:
\begin{enumerate}
    \item every point of $\mathbb{Z}_{2}^{n}$ lies on at least one line in $\mathcal{C}$,
    \item the underlying graph of $\mathcal{C}$ contains no even cycle, and
    \item no connected component of $\mathcal{C}$ is a tree.
\end{enumerate}
Then the X-ray transform restricted to $\mathcal{C}$ is injective; equivalently, $\mathcal{C}$ is admissible.
\end{theorem}

Moreover, we prove that within each connected component, at most one odd cycle may occur.

\begin{theorem}\label{thm:two-cycles}
Let $\mathcal{C}$ be a collection of lines in $\mathbb{Z}_{2}^{n}$, and suppose that a connected component of the underlying graph contains two distinct odd cycles. Then $\mathcal{C}$ is inadmissible.
\end{theorem}

Combining Theorems~\ref{thm:structure-overview} and \ref{thm:two-cycles}, we obtain a precise description of admissible complexes in $\mathbb{Z}_{2}^{4}$: each connected component consists of a single odd cycle with trees attached to its vertices. Since there are $16$ vertices in $\mathbb{Z}_{2}^{4}$, any odd cycle has length $3,5,7,9,11,13,$ or $15$.

Our second main result is an exact enumeration of admissible complexes in $\mathbb{Z}_{2}^{4}$ using this structural description and Cayley's formula for labeled trees.

\begin{theorem}\label{thm:main-count}
The total number of admissible line complexes in $\mathbb{Z}_{2}^{4}$ is
\[
984014621487058560 \;=\; 9.84014621487058560 \times 10^{17}.
\]
\end{theorem}

Finally, we describe and implement an algorithm that generalizes this counting method to $\mathbb{Z}_{2}^{n}$ for arbitrary $n$. As an illustration, we compute the number of admissible complexes in $\mathbb{Z}_{2}^{5}$.

\subsection{Relation to previous work}

The admissibility problem originated in the work of Gelfand and his collaborators \cite{gelfand1966}, in a continuous setting where the Radon transform integrates functions over families of affine subspaces in $\mathbb{R}^{n}$. Grinberg's thesis \cite{grinbergPHD} studied structural properties of line families, including their parametrization as manifolds via Pl\"{u}cker coordinates, and connected the admissibility problem to the geometry of symmetric spaces. Fundamental analytic and geometric aspects of the Radon transform and related transforms are developed in Helgason's monograph \cite{helgasonRadon} and in the tomography-focused exposition of Natterer \cite{nattererCT}.

In a different direction, Zentz \cite{zentz2009} analyzed admissible subsets of hyperplane spreads in finite affine spaces, determining minimal selections that preserve injectivity for hyperplane transforms. This work focuses on $k$-plane transforms in the finite setting and gives a classification of admissible hyperplane spreads in terms of combinatorial operations on parallel families. Our results are complementary: instead of hyperplanes, we treat the spatial X-ray transform on lines, and we develop a complete structural and enumerative description in the four-dimensional case, together with an algorithmic generalization to higher dimensions.

In the discrete case, Grinberg and others have investigated small dimensional examples such as $\mathbb{Z}_{2}^{3}$, both combinatorially and computationally; see \cite{Grinberg_2012} and the detailed catalog of inadmissible configurations in \cite{grinberg2021}. The scrapbook in particular highlights scenarios in which trees, even cycles, and more intricate obstructions destroy admissibility, as well as provides a blueprint for how local graph-theoretic features can be turned into global constraints on line complexes. 

Our work extends these ideas to dimension $4$ and beyond, replacing brute-force enumeration with a general structural framework and an algorithm that applies to $\mathbb{Z}_{2}^{n}$. The key significance of our results lies in the novel algorithm developed which solves the admissibility problem on $\mathbb{Z}_{2}^{n}$ for any number of dimensions $n$, as well as the implications for larger cases.


\subsection{Organization of the paper}

In Section~\ref{sec:preliminaries} we introduce the discrete Radon transform, line complexes, and the graph-theoretic language used throughout the paper. Section~\ref{sec:structure} proves the structural results on admissibility, including the obstructions given by omitted points, isolated trees, and even cycles, and the restriction to at most one odd cycle per component. In Section~\ref{sec:enumeration} we apply this structure to classify and enumerate admissible complexes in $\mathbb{Z}_{2}^{4}$, using Cayley's formula and a case division according to odd cycle lengths.

Section~\ref{sec:algorithm} describes an algorithm that generalizes the enumeration to $\mathbb{Z}_{2}^{n}$ and presents the results for $\mathbb{Z}_{2}^{5}$. In Section~\ref{sec:validation} we validate our counts using randomized sampling and rank computations of incidence matrices. Finally, Section~\ref{sec:discussion} discusses patterns in the prime factorizations of the counts, possible connections to spectral graph theory and singular value decomposition (SVD), and directions for future research.

\section{Preliminaries and notation}\label{sec:preliminaries}

\subsection{Finite vector spaces and lines}

We work over the two-element field $\mathbb{Z}_{2}$ throughout. For $n \geq 1$, the vector space $\mathbb{Z}_{2}^{n}$ consists of all $n$-tuples of elements of $\{0,1\}$, with coordinate-wise addition modulo $2$. We write $V_{n} = \mathbb{Z}_{2}^{n}$ for brevity, and denote $|V_{n}| = 2^{n}$.

\begin{figure}[h]
    \centering
    \includegraphics[scale = 0.12]{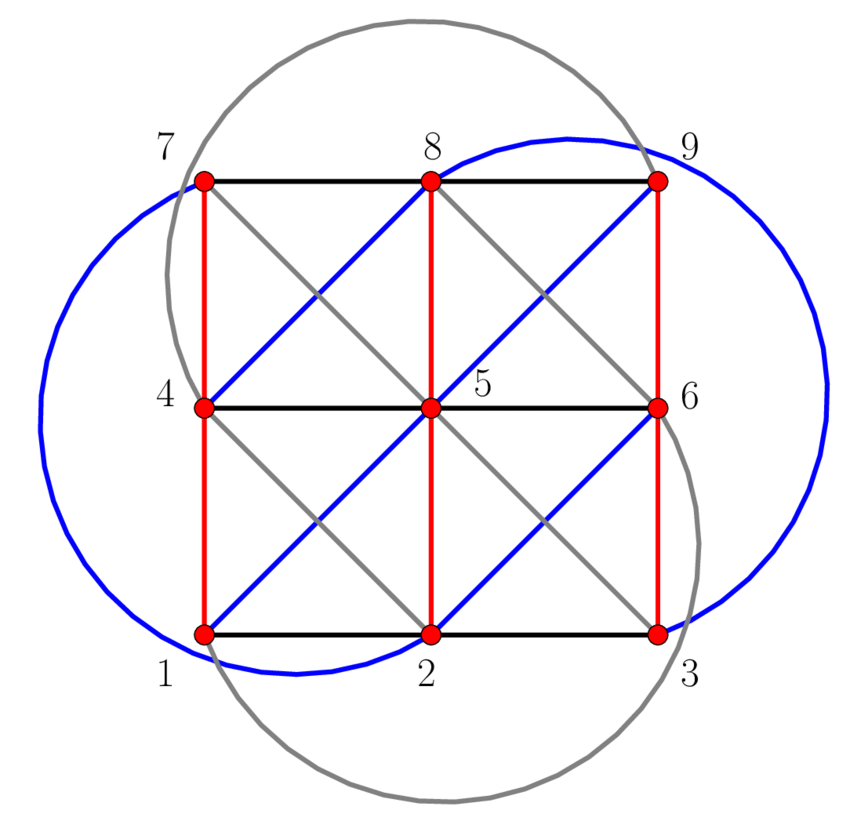}
    \caption{The affine plane of order $3$.}
    \label{fig:ctscan}
\end{figure}

\begin{definition}
An \emph{affine line} in $V_{n}$ is a subset of the form
\[
L = \{x, x+v\}
\]
where $x \in V_{n}$ and $v \in V_{n} \setminus \{0\}$. We identify lines that differ by a reparametrization, so each unordered pair of distinct points in $V_{n}$ determines a unique line.
\end{definition}

Thus the total number of lines in $V_{n}$ is $\binom{2^{n}}{2}$. We call any subset $\mathcal{C}$ of the set of all lines a \emph{line complex}.

Finite affine planes are referred to as $k$-planes for simplicity. In our space of interest, $\mathbb{Z}_{2}^{n}$, we integrate over ``2-planes'', or lines. Selections of lines in $\mathbb{Z}_{2}^{n}$ are the line complexes used to determine admissibility.

\subsection{The discrete X-ray transform}

We now define the discrete X-ray (line) transform on $V_{n}$.

\begin{definition}\label{def:radon}
Let $f \colon V_{n} \to \mathbb{C}$ be a function and let $L$ be an affine line in $V_{n}$. The \emph{discrete X-ray transform} of $f$ along lines is the function $R$ defined by
\[
R f(L) \;=\; \sum_{x \in L} f(x).
\]
If $\mathcal{C}$ is a set of lines in $V_{n}$, we denote by $R_{\mathcal{C}} f$ the collection $\{Rf(L)\}_{L \in \mathcal{C}}$.
\end{definition}

We can view $R_{\mathcal{C}}$ as a linear map from the vector space of functions $f \colon V_{n} \to \mathbb{C}$ to $\mathbb{C}^{|\mathcal{C}|}$. Choosing an ordering of the points and lines, we obtain a matrix representation of $R_{\mathcal{C}}$.

\begin{definition}
A line complex $\mathcal{C}$ in $V_{n}$ is \emph{admissible} if the map
\[
R_{\mathcal{C}} \colon \mathbb{C}^{V_{n}} \longrightarrow \mathbb{C}^{\mathcal{C}}
\]
is injective.
\end{definition}

Equivalently, $\mathcal{C}$ is admissible if the only function with $R_{\mathcal{C}} f = 0$ is the zero function.

\subsection{Graph-theoretic language}

The geometry of $V_{n}$ and its lines admits a convenient graph-theoretic interpretation, which mirrors the viewpoint taken in \cite{Grinberg_2012,grinberg2021}.

\begin{definition}
Given a line complex $\mathcal{C}$ in $V_{n}$, we define the associated graph $G(\mathcal{C})$ as follows:
\begin{itemize}
    \item the vertex set of $G(\mathcal{C})$ is $V_{n}$,
    \item two distinct vertices $x,y \in V_{n}$ are joined by an edge if and only if the line determined by $\{x,y\}$ lies in $\mathcal{C}$.
\end{itemize}
We use the terms \emph{vertex} and \emph{point} interchangeably.
\end{definition}

In this language, admissibility is a property of a spanning subgraph of the complete graph on $2^{n}$ vertices. We use standard graph-theoretic terminology.

\begin{definition}
A subgraph $T$ of $G(\mathcal{C})$ is a \emph{tree} if it is connected and acyclic. A \emph{cycle} is a simple closed path. A cycle is \emph{even} or \emph{odd} according to its length.
\end{definition}

We say that a vertex is \emph{omitted} by $\mathcal{C}$ if it has degree $0$ in $G(\mathcal{C})$, i.e., it lies on no line in $\mathcal{C}$.

\section{Structural obstructions to admissibility}\label{sec:structure}

In this section we identify simple graph-theoretic features of $G(\mathcal{C})$ that obstruct admissibility. These notions were already present in the study of $\mathbb{Z}_{2}^{3}$; we recall them and then extend the discussion to general $n$.

\subsection{Omitted points}

We begin with the simplest obstruction: a vertex not incident to any line of the complex.

\begin{lemma}\label{lem:omitted-point}
If a line complex $\mathcal{C}$ in $V_{n}$ omits a point, then $\mathcal{C}$ is inadmissible.
\end{lemma}

\begin{proof}
Suppose $p \in V_{n}$ lies on no line of $\mathcal{C}$. Define a function $f \colon V_{n} \to \mathbb{C}$ by $f(p) = 1$ and $f(x) = 0$ for all $x \neq p$. Then for every line $L \in \mathcal{C}$ we have $p \notin L$, so $Rf(L) = 0$. Hence $R_{\mathcal{C}} f = 0$ while $f \neq 0$. Thus $R_{\mathcal{C}}$ is not injective.
\end{proof}

\subsection{Isolated trees}

A second obstruction arises from tree components.

\begin{lemma}\label{lem:tree}
Let $\mathcal{C}$ be a line complex in $V_{n}$, and suppose that a connected component of $G(\mathcal{C})$ is a tree. Then $\mathcal{C}$ is inadmissible.
\end{lemma}

\begin{proof}
Let $T$ be a tree component with vertex set $V(T)$ and edge set $E(T)$. Since $T$ is a tree, $|E(T)| = |V(T)|-1$. Consider the restriction of the X-ray transform to the lines in $T$. The corresponding incidence matrix $M_{T}$ has one row for each edge and one column for each vertex, with entries in $\{0,1\}$ over $\mathbb{C}$. As a graph incidence matrix, $\mathrm{rank}(M_{T}) \leq |V(T)|-1$, since the sum of all columns is $0$.

Thus there exists a nonzero vector $f \in \mathbb{C}^{V(T)}$ in the kernel of $M_{T}$. Extend $f$ to all of $V_{n}$ by setting it to zero outside $V(T)$. Then $R_{\mathcal{C}}f$ vanishes on every line in $T$, and $f$ vanishes outside $T$, so $R_{\mathcal{C}}f = 0$ while $f \neq 0$. Hence $\mathcal{C}$ is inadmissible.
\end{proof}

\subsection{Even cycles}

The next obstruction comes from even cycles. This is already visible in $\mathbb{Z}_{2}^{3}$, where the $4$-cycle is the most intricate inadmissible configuration \cite{Grinberg_2012,grinberg2021}.

\begin{lemma}\label{lem:even-cycle}
If $G(\mathcal{C})$ contains an even cycle, then $\mathcal{C}$ is inadmissible.
\end{lemma}

\begin{proof}
Let $C$ be an even cycle in $G(\mathcal{C})$ with vertices $v_{1},v_{2},\dots,v_{2m}$ in cyclic order and edges $e_{1} = v_{1}v_{2}$, $e_{2} = v_{2}v_{3}$, $\dots$, $e_{2m} = v_{2m}v_{1}$. Define a function $f$ on these vertices by $f(v_{j}) = (-1)^{j}$ for $1 \leq j \leq 2m$, and set $f$ to be zero on all other vertices.

Each edge $e_{j}$ joins vertices with opposite labels, so
\[
Rf(e_{j}) = f(v_{j}) + f(v_{j+1}) = (-1)^{j} + (-1)^{j+1} = 0.
\]
Thus $Rf(L) = 0$ for every line $L$ in the cycle. Since $f$ is supported on the cycle and $\mathcal{C}$ contains all edges of $C$, the restriction of $R_{\mathcal{C}}$ to the cycle has a nontrivial kernel, and hence so does $R_{\mathcal{C}}$. Therefore $\mathcal{C}$ is inadmissible.
\end{proof}

\begin{figure}[h]%
    \centering
    \subfloat[]{{\includegraphics[width=5cm]{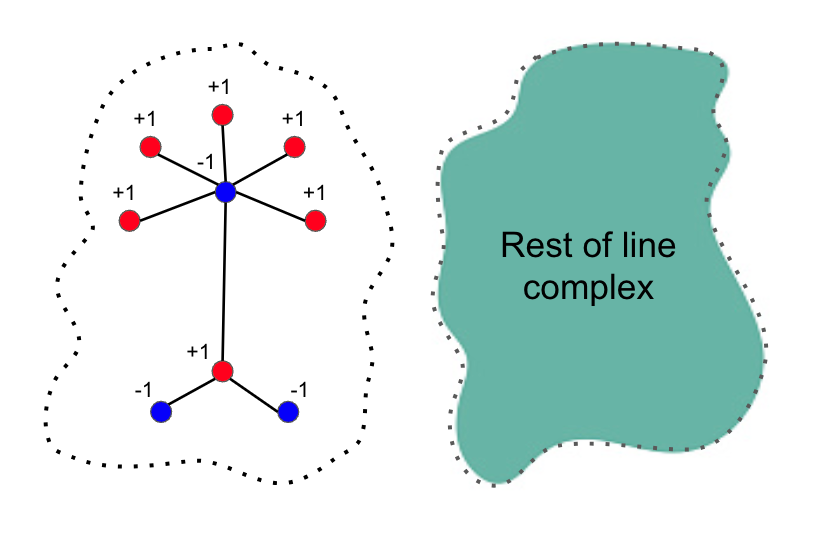} }}%
    \qquad
    \subfloat[]{{\includegraphics[width=3cm]{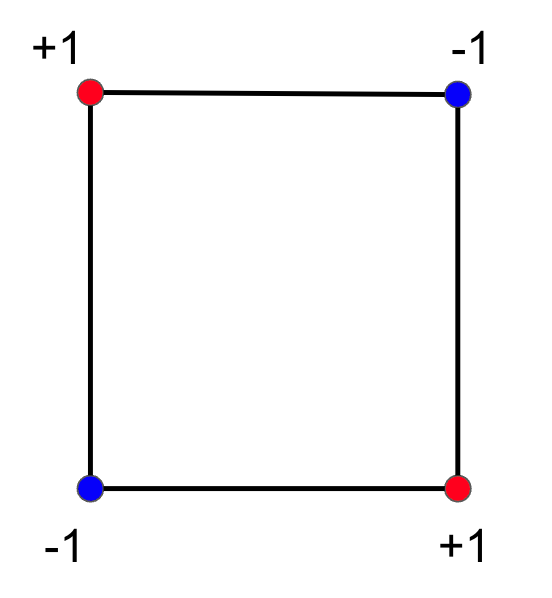} }}%
    \caption{Isolated trees (A) and even cycles (B) both fail to be admissible.}%
    \label{fig:example}%
\end{figure}

In $\mathbb{Z}_{2}^{3}$, the $4$-cycle is the largest possible even cycle and plays a special role: it generates a global linear dependency among line sums. Alternating vertex labels $\pm 1$ produce a function whose X-ray data vanishes on all four lines simultaneously, giving the most algebraically intricate inadmissible configuration in that space \cite{grinberg2021}.

\subsection{Odd cycles and self-inversion}

Odd cycles, in contrast, have a favorable property: the associated system of line sums is ``self-inverting.'' This observation underlies the proof of Theorem~\ref{thm:structure-overview}.

\begin{lemma}\label{lem:odd-cycle}
Let $C$ be an odd cycle with vertices $v_{1},\dots,v_{2m+1}$ and edges $e_{j} = v_{j}v_{j+1}$ (with indices taken modulo $2m+1$). Then the values $\{f(v_{j})\}$ can be recovered from the line sums $\{Rf(e_{j})\}$ by solving a nondegenerate linear system.
\end{lemma}

\begin{proof}
The system
\[
Rf(e_{j}) = f(v_{j}) + f(v_{j+1}), \quad j = 1,\dots,2m+1,
\]
is a linear system in the unknowns $f(v_{1}),\dots,f(v_{2m+1})$. The associated matrix is the incidence matrix of the cycle with entries in $\{0,1\}$. When the cycle length is odd, this matrix has full rank $2m+1$ over $\mathbb{C}$. One way to see this is to note that the only linear dependency among the edges of a cycle graph is the sum of all edges, which yields
\[
\sum_{j=1}^{2m+1} (f(v_{j}) + f(v_{j+1})) = 2 \sum_{j=1}^{2m+1} f(v_{j}),
\]
and in the odd-length case this does not vanish identically unless all $f(v_{j})$ vanish. It follows that the system is invertible.
\end{proof}

\subsection{Proof of the structural criterion}

We now prove Theorem~\ref{thm:structure-overview}, restated in a slightly more explicit form.

\begin{theorem}[Admissibility in $\mathbb{Z}_{2}^{n}$]\label{thm:admissibility-zn}
Let $\mathcal{C}$ be a line complex in $V_{n}$ such that:
\begin{enumerate}
    \item every vertex lies on at least one line in $\mathcal{C}$,
    \item $G(\mathcal{C})$ contains no even cycle, and
    \item no connected component of $G(\mathcal{C})$ is a tree.
\end{enumerate}
Then $\mathcal{C}$ is admissible.
\end{theorem}

\begin{proof}
Fix a connected component $H$ of $G(\mathcal{C})$. Since $H$ is not a tree and has no even cycle, it must contain at least one odd cycle. Let $C$ be a maximal odd cycle in $H$, and let $T$ denote the subgraph formed by the edges and vertices attached to $C$, excluding edges of $C$ itself. By construction, $T$ is a forest whose components are trees attached to vertices of $C$, because any additional cycle would necessarily introduce an even cycle or a second odd cycle within the component.

Consider a function $f$ supported on the vertices of $H$, and suppose that $R_{\mathcal{C}}f$ vanishes on all lines of $\mathcal{C}$ contained in $H$. We first restrict to the cycle $C$. By Lemma~\ref{lem:odd-cycle}, the line sums along the edges of $C$ determine all values of $f$ on the cycle. Since these line sums vanish, we conclude that $f$ vanishes on $C$.

Now consider a tree attached to a vertex $v$ of $C$. Every edge in the tree yields a linear equation of the form $f(u) + f(w) = 0$ for adjacent vertices $u,w$ in the tree. Starting from the boundary vertex $v$ where $f(v) = 0$, we can propagate along the tree and conclude that $f$ vanishes on all vertices in that tree. Repeating this argument for each attached tree shows that $f$ vanishes identically on $H$.

Since the argument applies to each connected component separately and every vertex of $V_{n}$ lies in some component, we obtain $f \equiv 0$ on all of $V_{n}$. Hence the kernel of $R_{\mathcal{C}}$ is trivial and $\mathcal{C}$ is admissible.
\end{proof}

We now prove the result restricting the number of odd cycles per connected component.

\begin{theorem}\label{thm:two-odd-cycles}
Let $\mathcal{C}$ be a line complex in $V_{n}$ whose underlying graph contains two odd cycles $A$ and $B$ in the same connected component. Then $\mathcal{C}$ is inadmissible.
\end{theorem}

\begin{proof}
Since $A$ and $B$ lie in the same connected component, there exists a path
\[
P = (p_{0},p_{1},\dots,p_{m})
\]
in $G(\mathcal{C})$ with $p_{0} \in A$ and $p_{m} \in B$, and edges $e_{1},\dots,e_{m}$ with $e_{j} = p_{j-1}p_{j}$. For each edge $e$ in $A$, $B$, or $P$, denote by $R_{e}$ the operation of summation over the corresponding line.

As in Lemma~\ref{lem:odd-cycle}, the line sums on $A$ (respectively $B$) determine the values of $f$ on that cycle. Consider the alternating sum
\[
S = R_{e_{1}} - R_{e_{2}} + R_{e_{3}} - \cdots \pm R_{e_{m}}.
\]
Each interior vertex of the path appears in exactly two edges of $P$ with opposite signs, so it cancels in $S$. Thus $S$ depends only on the values at the endpoints $p_{0}$ and $p_{m}$. More precisely, there is a sign $\sigma \in \{\pm 1\}$ such that
\[
S f = f(p_{0}) + \sigma f(p_{m}).
\]
But $f(p_{0})$ and $f(p_{m})$ can themselves be expressed in terms of line sums on $A$ and $B$, respectively. Hence $Sf$ is a linear combination of line sums over the edges of $A$ and $B$.

It follows that the collection of line sums consisting of those on $A$, those on $B$, and those on $P$ is linearly dependent. Therefore the corresponding rows in the incidence matrix of $R_{\mathcal{C}}$ are linearly dependent, and $R_{\mathcal{C}}$ cannot have full rank. Hence $\mathcal{C}$ is inadmissible.
\end{proof}

\begin{remark}
In particular, in an admissible complex every connected component contains exactly one odd cycle, with trees attached to its vertices. This cycle--tree decomposition is central to our enumeration in $\mathbb{Z}_{2}^{4}$ and mirrors, in four dimensions, the qualitative picture of admissible and inadmissible complexes illustrated in \cite{grinberg2021}.
\end{remark}

\begin{figure}[h]
    \centering
    \includegraphics[scale = 0.45]{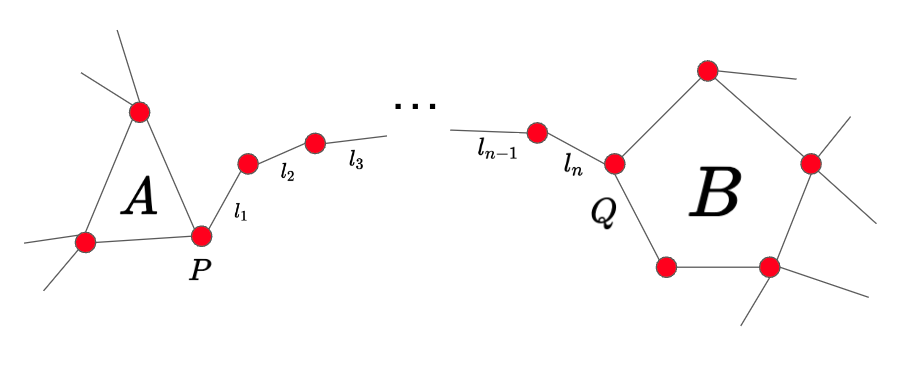}
    \caption{Visual representation of Theorem \ref{thm:two-odd-cycles}.}
    \label{fig:ctscan}
\end{figure}

\section{Classification and enumeration in $\mathbb{Z}_{2}^{4}$}\label{sec:enumeration}

We now specialize to the case $n=4$. The vector space $V_{4} = \mathbb{Z}_{2}^{4}$ has $16$ points and $\binom{16}{2} = 120$ lines. An admissible line complex must have exactly $16$ lines: one for each point in the sense of the X-ray transform restricted to a minimal invertible system, in analogy with the minimal admissible configurations analyzed in \cite{Grinberg_2012,grinberg2021}.

\subsection{Odd cycles in $\mathbb{Z}_{2}^{4}$}

As observed above, in an admissible complex $G(\mathcal{C})$ each connected component consists of a single odd cycle with attached trees. Since there are $16$ vertices in total, the possible cycle lengths are the odd integers between $3$ and $15$.

\begin{definition}
We say that an admissible complex in $\mathbb{Z}_{2}^{4}$ has \emph{cycle type}
\[
\{ \ell_{1},\dots,\ell_{r} \}
\]
if its connected components contain odd cycles of lengths $\ell_{1},\dots,\ell_{r}$, with $\sum_{j=1}^{r} \ell_{j} \leq 16$, and the remaining vertices are used to form trees attached to these cycles.
\end{definition}

The possible cycle types are listed in Table~\ref{tab:cycle-types}. Here we suppress the order of the $\ell_{j}$ and simply list all partitions of an odd integer $\leq 15$ into odd summands $\geq 3$.

\begin{table}[h]
\centering
\renewcommand{\arraystretch}{1.15}
\begin{tabular}{|c|c|c|c|c|}
\hline
\textbf{1 cycle} & \textbf{2 cycles}  & \textbf{3 cycles}  & \textbf{4 cycles}      & \textbf{5 cycles}          \\ \hline
$\{3\}$   & $\{3,3\}$   & $\{3,3,3\}$ & $\{3,3,3,3\}$ & $\{3,3,3,3,3\}$ \\
$\{5\}$   & $\{3,5\}$  & $\{3,3,5\}$ & $\{3,3,3,5\}$ &                   \\
$\{7\}$   & $\{3,7\}$  & $\{3,3,7\}$ & $\{3,3,3,7\}$ &                   \\
$\{9\}$   & $\{3,9\}$  & $\{3,3,9\}$ & $\{3,3,5,5\}$ &                   \\
$\{11\}$  & $\{3,11\}$ & $\{3,5,5\}$ &                &                   \\
$\{13\}$  & $\{3,13\}$ & $\{3,5,7\}$ &                &                   \\
$\{15\}$  & $\{5,5\}$  & $\{5,5,5\}$ &                &                   \\
          & $\{5,7\}$  &             &                &                   \\
          & $\{5,9\}$  &             &                &                   \\
          & $\{5,11\}$ &             &                &                   \\
          & $\{7,7\}$  &             &                &                   \\
          & $\{7,9\}$  &             &                &                   \\
\hline
\end{tabular}
\caption{Possible combinations of disjoint odd cycles in $\mathbb{Z}_{2}^{4}$.}
\label{tab:cycle-types}
\end{table}

For each cycle type, we must:
\begin{enumerate}
    \item choose the vertices that belong to the cycles,
    \item count the number of ways to form the cycles on those vertices,
    \item distribute the remaining vertices among the cycle vertices, and
    \item count the number of trees that can be formed on the resulting vertex sets.
\end{enumerate}

\subsection{Cayley's formula and cycles}

Our enumeration uses Cayley's classical formula.

\begin{theorem}[Cayley's formula {\cite{cayley_original,shukla2017,proofsFromBook}}]\label{thm:cayley}
There are exactly $n^{n-2}$ labeled trees on $n$ vertices.
\end{theorem}

We also require the number of cycles on a given vertex set.

\begin{proposition}\label{prop:cycles}
Let $n \geq 3$. The number of distinct cycles that can be formed on a set of $n$ labeled vertices is
\[
\frac{(n-1)!}{2}.
\]
\end{proposition}

\begin{proof}
There are $n!$ permutations of the $n$ vertices. Each cycle can be represented by $2n$ of these permutations, since we can start at any of the $n$ vertices and traverse the cycle in either direction. Thus the number of distinct cycles is $n!/(2n) = (n-1)!/2$.
\end{proof}

\subsection{Example: the $\{3,3\}$ case}

We illustrate the counting procedure in the case of two disjoint $3$-cycles. There are $6$ vertices in the cycles and $10$ remaining vertices to be attached as trees.

We first choose the $6$ cycle vertices:
\[
\binom{16}{3} \binom{13}{3}
\]
ways, and by Proposition~\ref{prop:cycles} there is exactly one cycle on each set of three vertices, since $(3-1)!/2 = 1$.

Next we must distribute the $10$ remaining vertices among the $6$ cycle vertices. Let $A_{1},\dots,A_{6}$ denote the number of additional vertices attached to each cycle vertex. Then
\[
A_{1} + \cdots + A_{6} = 10,
\]
and each ordered $6$-tuple $(A_{1},\dots,A_{6})$ with nonnegative integer entries satisfying this equation corresponds to a way of assigning the extra vertices.

The possible distributions $(A_{1},\dots,A_{6})$ are all possible assignments of each $0\le A_i\le 10$ such that $\sum_{i=1}^6 A_i=10$. For each such distribution, we must:
\begin{itemize}
    \item choose which of the $10$ extra vertices are attached to which cycle vertex (up to relabeling), and
    \item count the number of trees on each group of vertices using Cayley's formula.
\end{itemize}

For each distribution we compute a multiplier accounting for symmetry (repeated entries among the $A_{i}$). We then apply Cayley's formula to count labeled trees on each group consisting of a cycle vertex together with its attached vertices, and we multiply over all six groups.

Summing over all distributions, and then multiplying by the number of ways to choose the cycle vertices, we find that the total number of admissible complexes of cycle type $\{3,3\}$ is:
\[
66036668364226560 = 2^{42} \cdot 3 \cdot 5 \cdot 7 \cdot 11 \cdot 13.
\]

\subsection{Enumeration for all cycle types}

We carry out an analogous computation for each cycle type listed in Table~\ref{tab:cycle-types}. For one-cycle types, the count divides into:
\begin{itemize}
    \item choosing the $\ell$ cycle vertices from $16$,
    \item counting the number of cycles on them (Proposition~\ref{prop:cycles}), and
    \item distributing the remaining $16-\ell$ vertices and forming trees using Theorem~\ref{thm:cayley}.
\end{itemize}

For multiple-cycle types, there is an additional combinatorial factor accounting for the choice of disjoint vertex sets for each cycle.

Table~\ref{tab:cases-1cycle} lists the contribution of each one-cycle case and the corresponding prime factorization.

\noindent \begin{table}[h]
\centering
\renewcommand{\arraystretch}{1.25}
\begin{tabular}{|l|l|l|l|l|}
\hline
Cycle Type & Initial Count & Multiplier & Final Total & Prime Factorization \\ \hline
$\{3\}$        & $844424930131968$        & $\binom{16}{3}$        & $472877960873902080$
        & $2^{52} \cdot 3 \cdot 5 \cdot 7$            \\ 
$\{5\}$        & $5497558138880$        & $\binom{16}{5} \cdot \frac{4!}{2}$         & $288160007407534080$       & $2^{46} \cdot 3^2 \cdot 5 \cdot 7 \cdot 13$         \\ 
$\{7\}$        & $30064771072$       & $\binom{16}{7} \cdot \frac{6!}{2}$        & $123818753182924800$        & $2^{39} \cdot 3^2 \cdot 5^2 \cdot 7 \cdot 11 \cdot 13$         \\ 
$\{9\}$        & $150994944$        & $\binom{16}{9} \cdot \frac{8!}{2}$        & $34824024332697600$        & $2^{34} \cdot 3^4 \cdot 5^2 \cdot 7 \cdot 11 \cdot 13$            \\ 
$\{11\}$        & $720896$        & $\binom{16}{11} \cdot \frac{10!}{2}$        & $5713316492083200$        & $2^{27} \cdot 3^5 \cdot 5^2 \cdot 7^2 \cdot 11 \cdot 13$              \\ 
$\{13\}$        & $3328$        & $\binom{16}{13} \cdot \frac{12!}{2}$        & $446352850944000$        & $2^{21} \cdot 3^5 \cdot 5^3 \cdot 7^2 \cdot 11 \cdot 13$               \\ 
$\{15\}$        & $15$        & $\binom{16}{15} \cdot \frac{14!}{2}$        & $10461394944000$        & $2^{14} \cdot 3^6 \cdot 5^3 \cdot 7^2 \cdot 11 \cdot 13$              \\ 
\hline
\end{tabular}
\caption{Total admissible complexes for one-cycle types in $\mathbb{Z}_{2}^{4}$.}
\label{tab:cases-1cycle}
\end{table}

Table~\ref{tab:cases-2plus} begins the list of multi-cycle types, starting with the $\{3,3\}$ case discussed above. The remaining cases are similarly computed.

\begin{table}[h]
\centering
\renewcommand{\arraystretch}{1.25}
\begin{tabular}{|l|l|l|l|l|}
\hline
Cycle Type & Initial Count & Multiplier & Final Total & Prime Factorization \\ \hline
$\{3,3\}$        & $412316860416$        &  $\binom{16}{3}\binom{13}{3}$       & $66036668364226560$        & $2^{42} \cdot 3 \cdot 5 \cdot 7 \cdot 11 \cdot 13$               \\  
$\vdots$        & $\vdots$        & $\vdots$        & $\vdots$        & $\vdots$               \\ 
\hline
\end{tabular}
\caption{Sample multi-cycle cases in $\mathbb{Z}_{2}^{4}$; the omitted rows correspond to other cycle types in Table~\ref{tab:cycle-types}.}
\label{tab:cases-2plus}
\end{table}

Summing over all cycle types, we obtain Theorem~\ref{thm:main-count}.

\begin{proof}[Proof of Theorem~\ref{thm:main-count}]
For each cycle type listed in Table~\ref{tab:cycle-types} we compute the number of admissible complexes by the procedure outlined above: choose cycle vertices, form cycles, distribute remaining vertices, count trees via Cayley's formula, and account for symmetry. Summing the resulting counts over all cycle types yields
\[
984014621487058560\approx 9.840\times 10^{17}
\]
admissible complexes in $\mathbb{Z}_{2}^{4}$. The combinatorial features of these counts, such as their prime factorizations in Tables~\ref{tab:cases-1cycle} and \ref{tab:cases-2plus}, suggest that a more direct combinatorial parametrization of admissible complexes may exist.
\end{proof}

\section{An algorithm for $\mathbb{Z}_{2}^{n}$}\label{sec:algorithm}

The enumeration method described in Section~\ref{sec:enumeration} extends naturally to $\mathbb{Z}_{2}^{n}$ for general $n$. The key observation is that admissible complexes are determined by:
\begin{itemize}
    \item a partition of an odd integer $\leq 2^{n}$ vertices into odd cycle lengths, and
    \item the trees attached to the cycle vertices.
\end{itemize}

\subsection{Partitioning into odd cycle lengths}

Let $N = 2^{n}$. Any admissible complex decomposes as a disjoint union of connected components, each containing one odd cycle. If the cycle lengths are $\ell_{1},\dots,\ell_{r}$, then
\[
\ell_{1} + \cdots + \ell_{r} \leq N,
\]
and each $\ell_{j}$ is an odd integer at least $3$. Thus the problem of enumerating cycle types reduces to finding all partitions of odd integers $\leq N$ into odd parts $\geq 3$.

We implement this using a simple recursive procedure. Let $P(j)$ denote the set of partitions of $j$ into odd parts, and build up $P(j)$ for $j=0,1,\dots,N$ by adding odd parts one at a time. The resulting list $P(N)$ encodes all possible cycle types.

\subsection{Counting trees via Cayley's formula}

Once a cycle type $(\ell_{1},\dots,\ell_{r})$ is fixed, we proceed as in Section~\ref{sec:enumeration}. We:
\begin{enumerate}
    \item choose disjoint subsets of $V_{n}$ of sizes $\ell_{1},\dots,\ell_{r}$ to serve as cycle vertices,
    \item count the number of cycles on each subset using Proposition~\ref{prop:cycles},
    \item distribute the remaining $N - \sum_{j}\ell_{j}$ vertices among the cycle vertices, and
    \item count the number of labeled trees that can be formed on each group using Theorem~\ref{thm:cayley}.
\end{enumerate}

The only difference from the $n=4$ case is that the total number of vertices $N$ is now larger, leading to more partitions and more distributions per cycle type. Nonetheless, for moderate values of $n$ (e.g., $n \leq 5$) the procedure is computationally tractable. The entire algorithm is presented in Algorithm \ref{alg:cap}.

\begin{algorithm}[h]
\caption{Partitioning values $\le n$ into odd sums}\label{alg:cap}
\begin{algorithmic}
\Require $n \geq 0$ and $n$ even
\State arr $\gets [\hspace{5pt}]$
\For{\texttt{i in range(n+1)}}
    \State \texttt{arr $\gets [\hspace{5pt}]$}
\EndFor
\State arr[0] $\gets [\hspace{2pt}[\hspace{4pt}]\hspace{2pt}]$
\For{\texttt{j in range(1, n+1)}}
    \For{\texttt{k in range(1, j+1, 2)}}
        \For{\texttt{partition in arr[j-k]}}
            \State \texttt{arr[j] $\gets $partition+[j]}
        \EndFor
    \EndFor
\EndFor
\State \textbf{return} arr[n] \Comment{We return the entire set of partitions.}
\end{algorithmic}
\end{algorithm}

\subsection{Results for $\mathbb{Z}_{2}^{5}$}

We implemented this algorithm to compute the number of admissible line complexes in $\mathbb{Z}_{2}^{5}$. The space $V_{5} = \mathbb{Z}_{2}^{5}$ has $32$ points and $\binom{32}{2} = 496$ lines. The total number of possible line complexes of size $32$ is
\[
\binom{496}{32} \approx 2.463 \times 10^{50},
\]
so brute-force enumeration is infeasible. Our algorithm, however, runs in well under a minute for $n=5$. The results for $n=3,4,5$ are summarized in Table~\ref{tab:summary-n}.

\begin{table}[h]
\centering
\renewcommand{\arraystretch}{1.25}
\begin{tabular}{|c|c|c|c|c|c|}
\hline
Space & Points & Lines & Total Complexes & Admissible & Percent Admissible \\ \hline
$\mathbb{Z}_{2}^{3}$ & $8$  & $28$   & $3.108 \cdot 10^{6}$  & $9.3744 \cdot 10^{5}$ & $30\%$ \\
$\mathbb{Z}_{2}^{4}$ & $16$ & $120$  & $3.104 \cdot 10^{19}$ & $9.840 \cdot 10^{17}$ & $3\%$  \\  
$\mathbb{Z}_{2}^{5}$ & $32$ & $496$  & $2.463 \cdot 10^{50}$ & $6.817 \cdot 10^{46}$ & $0.03\%$ \\ \hline
\end{tabular}
\caption{Summary of admissible complexes in $\mathbb{Z}_{2}^{n}$ for $n=3,4,5$. The values for $n=3$ and $n=4$ agree with independent computations and previous work, while the $n=5$ values are obtained from our algorithm.}
\label{tab:summary-n}
\end{table}

As $n$ increases, the proportion of admissible complexes drops rapidly. This reflects the increasing number of ways to introduce inadmissible obstructions such as even cycles and tree components, as already suggested in the three-dimensional scrapbook \cite{grinberg2021}.

\section{Validation via randomized sampling}\label{sec:validation}

To validate our counts for $\mathbb{Z}_{2}^{4}$, we performed randomized sampling experiments. We generated $10^{6}$ random line complexes of size $16$ by choosing $16$ lines uniformly at random from the set of all $120$ lines in $V_{4}$, and we tested admissibility by computing the rank of the corresponding incidence matrices over $\mathbb{Q}$.

Over these $10^{6}$ random complexes, the proportion of admissible complexes was approximately 0.031505.\\
\\
Since there are $\binom{120}{16} \approx 3.104406 \times 10^{19}$ possible complexes, this suggests an expected number of admissible complexes
\[
(0.031505) \times (3.104406 \times 10^{19}) \approx 9.78043 \times 10^{17},
\]
which is consistent with our exact count
\[
9.84014621487058560 \times 10^{17}
\]
from Theorem~\ref{thm:main-count}. The small discrepancy is well within the expected sampling error and is further explained by the fact that complexes with multiple cycles are much rarer than single-cycle complexes.

\begin{remark}
An alternative validation method is to study the singular value decomposition (SVD) of the incidence matrices. Admissible complexes correspond to matrices with full column rank, and the smallest singular value is bounded away from zero. Inadmissible complexes exhibit singular values that are exactly zero. This spectral perspective may be useful for studying stability properties of the inverse problem and parallels the conditioning questions that arise in continuous tomography \cite{nattererCT}.
\end{remark}

\section{Discussion and future directions}\label{sec:discussion}

Our results reveal several intriguing patterns and suggest multiple directions for future work. Firstly, the prime factorizations appearing in Tables~\ref{tab:cases-1cycle} and \ref{tab:cases-2plus} display nontrivial regularities. The exponents of primes such as $2,3,5,7$, and $11$ vary systematically with the cycle type, suggesting that a more direct combinatorial parametrization of admissible complexes may exist. One possibility is to encode admissible complexes via generalized Prüfer sequences or other combinatorial objects that naturally produce powers of these primes, as in classical treatments of Cayley's theorem \cite{shukla2017,proofsFromBook}.

A long-term goal is to derive a closed-form expression, or at least an asymptotic formula, for the number of admissible complexes in $\mathbb{Z}_{2}^{n}$ as $n \to \infty$.

While our work focuses on the two-element field, it is natural to consider $\mathbb{Z}_{p}^{n}$ for odd primes $p$ or more general finite fields $\mathbb{F}_{q}$. In these settings, lines contain more than two points, and the X-ray transform sums over larger affine subspaces. The graph-theoretic picture would need to be replaced by a hypergraph model, and even/odd cycle obstructions would have to be generalized to higher-dimensional analogues. Analytic results for finite-field Radon transforms, such as those in \cite{koh2012sharp}, suggest that rich structure should persist in this more general setting.

Understanding how our cycle--tree decomposition extends to these settings could provide new insights into admissibility for higher-dimensional Radon transforms, both discrete and continuous.

\subsection{Connections to imaging and inverse problems}

Admissible complexes model sampling patterns for tomographic imaging: each line corresponds to a measurement of the sum of densities along a ray. In computed tomography (CT) and other imaging modalities, physical or technological constraints often prevent the collection of all possible line integrals, and one must instead choose an efficient and informative subset. The mathematical theory of such transforms has been developed in detail in \cite{helgasonRadon,nattererCT}.

Our structural criteria and enumeration results suggest that admissible sampling patterns in high-dimensional discrete models are relatively rare but can be characterized explicitly. Extending this viewpoint to continuous models, or to discretizations of practical imaging geometries, may help optimize scan designs and reduce redundancy. Figure 5 shows an explicit real-world example of admissible sampling in a continuous case.

\begin{figure}[h]
    \centering
    \includegraphics[scale = 0.15]{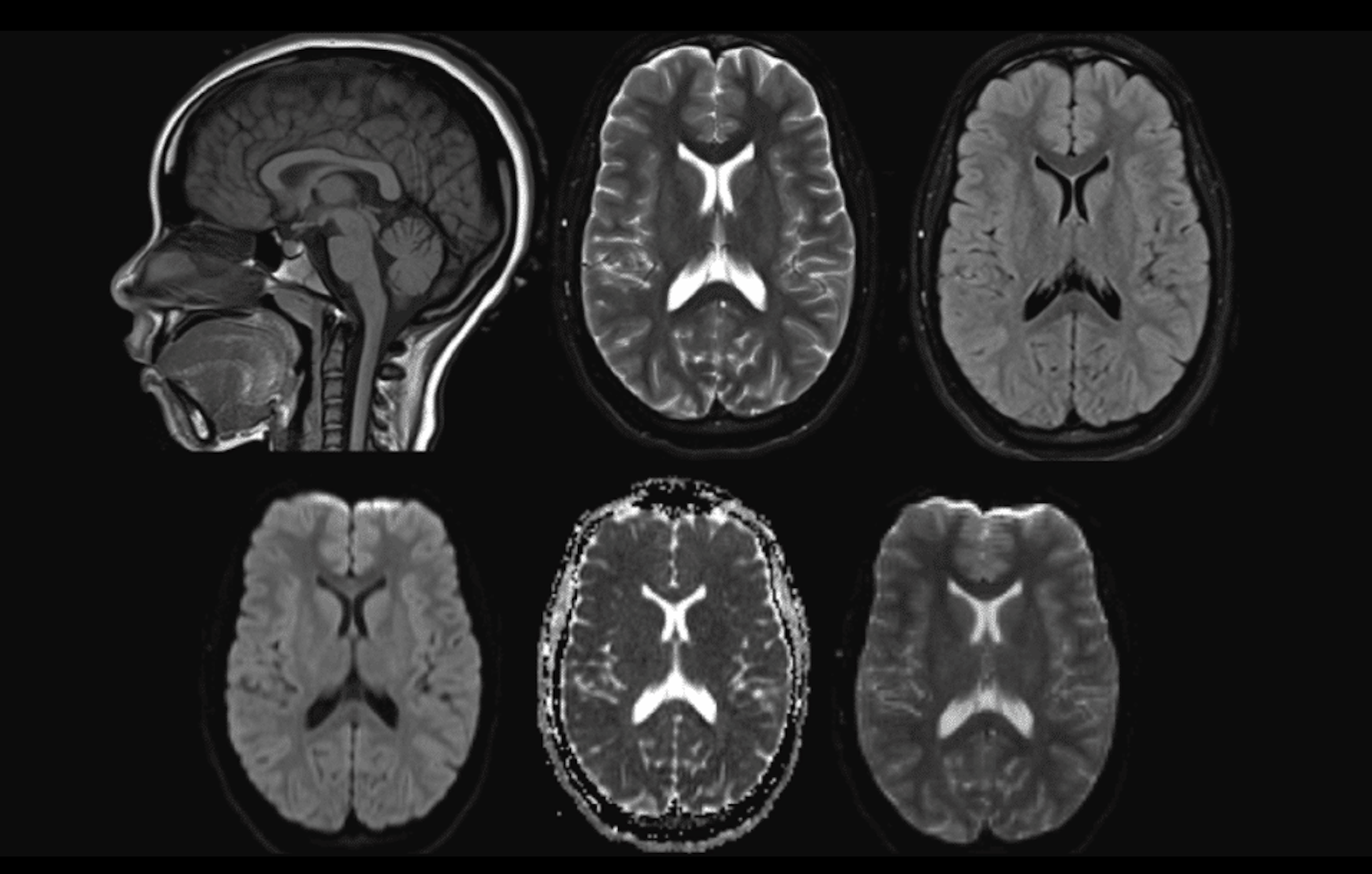}
    \caption{Planar cross-sections of the brain can be used to recover original tissue densities, due to the invertibility of the X-ray transform \cite{nattererCT}.}
    \label{fig:ctscan}
\end{figure}

\subsection{Spectral and algebraic perspectives}

Finally, the incidence matrices associated with admissible complexes form an interesting family of full-rank matrices with strong combinatorial structure. Studying their eigenvalues and singular values may shed light on the stability of the inversion process and the conditioning of the associated linear systems. It would be particularly interesting to connect these spectra to well-studied families of graphs or to representation-theoretic structures arising from the symmetry group of $V_{n}$, in the spirit of the group-theoretic viewpoint on Radon transforms emphasized in \cite{helgasonRadon}.

\section*{Acknowledgments} 

The first author would like to thank her mentor, Dr.\ Eric Grinberg, for his guidance throughout this project. She is grateful to Peter Gaydarov for valuable discussions about the admissibility problem and help in preparing this manuscript. She also thanks the Research Science Institute (RSI), the Center for Excellence in Education (CEE), and the Massachusetts Institute of Technology (MIT) for their support and for providing an environment in which this research could be conducted.

\bibliographystyle{amsplain}
\bibliography{biblio}

\end{document}